\RequirePackage{ifpdf}
\ifpdf % We are running pdfTeX in pdf mode
\documentclass[pdftex]{sigma}
\else
\documentclass{sigma}
\fi

\usepackage{enumerate}

\numberwithin{equation}{section}

\newcommand{\bp}{p}
\newcommand{\bu}{u}
\newcommand{\bv}{v}
\newcommand{\bH}{{\mathbf H}}
\newcommand{\bV}{{\mathbf V}}
\newcommand{\bE}{{\mathbf E}}

\newcommand{\bC}{\mathbb{C}}
\newcommand{\bN}{\mathbb{N}}
\newcommand{\bR}{\mathbb{R}}
\newcommand{\bZ}{\mathbb{Z}}

\newcommand{\Ht}{{\mathrm{ht}}}

\newcommand{\id}{{\mathrm{id}}}
\newcommand{\odd}{{\mathrm{odd}}}

\newcommand{\sgn}{{\mathrm{sgn}}}
\newcommand{\Tab}{{\mathrm{Tab}}}
\newcommand{\HVTab}{{\mathrm{Tab}_{\scriptscriptstyle \mathrm HV}}}
\newcommand{\HVP}{P_{\scriptscriptstyle \mathrm HV}}

\newcommand{\Pos}{{\mathrm{Pos}}}

\newcommand{\cA}{\mathcal{A}}

\newcommand{\cH}{\mathcal{H}}
\newcommand{\cLU}{\mathcal{LU}}
\newcommand{\cP}{\mathcal{P}}

\newcommand{\cTh}{\mathcal{T}_{\mathrm h}}
\newcommand{\cU}{\mathcal{U}}
\newcommand{\cW}{\mathcal{W}}
\newcommand{\cZ}{\mathcal{Z}}

\newcommand{\fg}{\mathfrak{g}}
\newcommand{\hg}{\hat{\fg}}
\newcommand{\fP}{\mathfrak{P}}
\newcommand{\fS}{\mathfrak{S}}

\newcommand{\Uqhg}{U_q(\hg )}

\newcommand{\pathNakanishi}[3]{#1 \overset{#2}{\to} #3}
\newcommand{\I}{\mathrm{I}}
\newcommand{\II}{\mathrm{II}}
\newtheorem{thm}{Theorem}[section]
\newtheorem{lem}[thm]{Lemma}
\newtheorem{prop}[thm]{Proposition}
\theoremstyle{definition}
\newtheorem{defn}[thm]{Definition}
\newtheorem{exmp}[thm]{Example}
\newtheorem{rem}[thm]{Remark}
\begin{document}

\allowdisplaybreaks

\renewcommand{\PaperNumber}{078}

\FirstPageHeading

\ShortArticleName{Jacobi--Trudi Determinant of $C_n$}

\ArticleName{Paths and Tableaux Descriptions\\ of
Jacobi--Trudi Determinant Associated\\ with
Quantum Af\/f\/ine Algebra of Type $\boldsymbol{C_n}$}

\Author{Wakako NAKAI and  Tomoki NAKANISHI}

\AuthorNameForHeading{W. Nakai and T. Nakanishi}

\Address{Graduate School of Mathematics,
Nagoya University, Nagoya 464-8602, Japan}

\Email{\href{mailto:m99013c@math.nagoya-u.ac.jp}{m99013c@math.nagoya-u.ac.jp}, \href{mailto:nakanisi@math.nagoya-u.ac.jp}{nakanisi@math.nagoya-u.ac.jp}}

\ArticleDates{Received May 03, 2007, in f\/inal form July 04, 2007; Published online July 18, 2007}

\Abstract{We study the Jacobi--Trudi-type determinant which is conjectured to be
the $q$-character of a certain, in many cases irreducible, f\/inite-dimensional
representation of the quantum af\/f\/ine algebra of type $C_n$.
Like the $D_n$ case studied by the authors recently,
applying the Gessel--Viennot path method with
 an additional involution
and a deformation of paths, we obtain an expression by a positive sum over
a set of tuples of paths,
which is naturally translated into the one over a set of tableaux
on a skew diagram.}


\Classification{17B37; 05E15}

%%%%%%%%%%%%%%%%%%%%%
\section{Introduction}
Let $\fg$ be a simple Lie algebra over $\bC$, and let
$\hg$ be the corresponding untwisted af\/f\/ine Lie algebra.
Let $\Uqhg$ be the quantum af\/f\/ine algebra, namely,
the quantized universal enveloping algebra
of~$\hg$~\cite{D1,J}.
In order to investigate the
f\/inite-dimensional representations of~$\Uqhg$~\cite{D2,CP},
the {\it  $q$-character} of $\Uqhg$ was introduced
and studied in \cite{FR,FM}.
So far, the explicit description of
$\chi_q(V)$ is available only for
a limited type of representations (e.g., the fundamental representations
\cite{FM,CM}),
and the description  for
general $V$ is an open problem.
There are also results  by \cite{N1,N2, H1, H2}, some of which
will be mentioned below.

This is a continuation of our previous work \cite{NN1,NN2}.
Firstly, let us brief\/ly summarize some background.
In \cite{NN1}, for $\fg$ of types $A_n$, $B_n$, $C_n$, and $D_n$,
we conjecture
that the $q$-characters
of a~certain family of, in many cases irreducible, f\/inite-dimensional
representations  of $\Uqhg$ are given by
the determinant form $\chi_{\lambda/\mu,a}$,
where $\lambda/\mu$ is a skew diagram and $a$ is a complex
parameter.
We call $\chi_{\lambda/\mu,a}$ the Jacobi--Trudi determinant.
For $A_n$ and $B_n$, this is a reinterpretation
of the conjecture for the spectra
of the transfer matrices of the vertex models
associated with the corresponding representations \cite{BR,KOS}.
When a diagram $\lambda/\mu$ is a rectangle of depth $d$
with $d\leq n$ for $A_n$, $d\leq n-1$ for $B_n$ and $C_n$,
and $d\leq n-2$ for $D_n$,
the corresponding representation $V$ is a so-called {\it Kirillov--Reshetikhin (KR) module}
\cite{KR, KN}.
In this case, the conjecture is known to be true from the facts:
(i) The determinants $\chi_{\lambda/\mu,a}$ corresponding to the KR modules
consist of  a part of the unique solution of the set of  polynomial relations called the
{\it $T$-system}
 (\cite[equation~(2.17)]{KNS} for $A_n$,
 \cite[Theorem 5.1]{KOS} for $B_n$,
  \cite[Theorem 3.1]{KNH} for $C_n$,
  \cite[Theorem 3.1]{TK} for $D_n$).
(ii) The $q$-characters of the KR modules
also consists of  the solution of the $T$-system (\cite[Theorem 1.1]{N2} for $A_n$ and $D_n$,
 \cite[Theorem 3.4]{H1} for any $X_n$).

The characters of an irreducible representations of $sl_n$,
or   the Schur polynomials,
are described  by the semistandard tableaux on Young diagrams.
As is well known, several objects of
representation theory and combinatorics
are related via  the semistandard tableaux (e.g.\  \cite{S}).
Since $\chi_{\lambda/\mu,a}$ is expected to be the characters of
representations of $\Uqhg$,
it is natural to ask if there are such tableaux for
$\chi_{\lambda/\mu,a}$.
This is the question studied in \cite{NN1,NN2}.
To do it, the paths method by  \cite{GV}  (see also \cite{S}),
which derives the semistandard tableaux purely combinatorially
from the determinant expression of the Schur polynomials,
seems suited to our situation.

Now let us brief\/ly summarize the results in \cite{NN1,NN2}.
For $A_n$ and $B_n$ \cite{NN1},
the application of the standard paths method
by \cite{GV}   immediately reproduces
the known tableaux descriptions of $\chi_{\lambda/\mu,a}$
by \cite{BR,KOS}.
Here,  by {\it tableaux description\/} we mean
an expression of
$\chi_{\lambda/\mu,a}$ by a positive sum
over a certain set of tableaux on $\lambda/\mu$.
In contrast, for $C_n$ \cite{NN1},
the situation is not as simple as the former cases,
and we obtain a tableaux description of $\chi_{\lambda/\mu,a}$
only for some special cases (i.e.,
a skew diagram $\lambda/\mu$ of at most two columns or
of at most three rows).
In \cite{NN2}, we consider the same problem for $D_n$,
where the situation is  similar to $C_n$.
By extending the idea of \cite{NN1} to full generality,
a tableaux description of $\chi_{\lambda/\mu,a}$
for a general skew diagram $\lambda/\mu$
is obtained.

In this paper, we come back to $C_n$ again,
and, using the method in \cite{NN2}, we now obtain
a~tableaux description of $\chi_{\lambda/\mu,a}$
for a general skew diagram $\lambda/\mu$.
Even though the formulation is quite parallel
to $D_n$,
the $C_n$ case has slight technical
complications in two aspects.
First, the def\/initions of the duals
of lower and upper paths in \eqref{eq:alpha-dual}
and \eqref{eq:beta-dual} are more delicate.
Secondly, the translation of the extra rule into
tableau language is less straightforward.
In retrospect, these are the reasons
why we could derive
the tableaux description in a general case f\/irst for $D_n$.
We expect that the method is also applicable
to the twisted quantum af\/f\/ine algebras of classical type \cite{T},
and hopefully, even to the quantum af\/f\/ine algebras
of exceptional type as well.
We also expect that our tableaux description
(and the corresponding one in the path picture)
is useful to study the $q$-characters and the crystal bases
\cite{KaN,OSS} of those representations.
One of the advantage of our tableaux (or paths) is that
it naturally resolves the multiplicity of the weight
polynomials of $q$-characters so that
the proposed algorithm by \cite{FM} to create
the $q$-characters could be more naturally realized
on the space of tableaux (or paths) for those representations
having multiplicities.

The organization of the paper is as follows.
In Section \ref{sec:JT-det} we review the result of \cite{NN1}.
Namely, we def\/ine the Jacobi--Trudi determinant $\chi_{\lambda/\mu,a}$
for $C_n$; then,  following the standard method by~\cite{GV},
we introduce lattice paths and express  $\chi_{\lambda/\mu, a}$
as a sum
over a certain set of tuples of paths.
In Section~\ref{sec:positive-sum},
we apply the method of~\cite{NN2},
with adequate modif\/ications of the basic notions for $C_n$,
and obtain an expression of $\chi_{\lambda/\mu, a}$
by a positive sum
over a set of tuples of paths (Theorem~\ref{thm:tilde-positive-sum}).
In Section~\ref{sec:tableau-description},
we  translate the last expression into the tableaux description
whose tableaux are determined by the
horizontal, vertical, and extra rules
(Theorems \ref{thm:tableau-description}
and~\ref{thm:configuration}).

This paper is written as the companion to~\cite{NN2}.
As we mentioned, once the basic objects are   adequately set,
the derivation of paths and tableaux descriptions is quite parallel to $D_n$.
In particular, the proofs of the core propositions in Section 3 are the same
almost word for word.
For such propositions, instead of repeating the proofs,
we only refer the corresponding propositions in \cite{NN2}
so that readers could focus on what are specif\/ic to the $C_n$ case.

\section[Jacobi-Trudi determinant and paths]{Jacobi--Trudi determinant and paths}\label{sec:JT-det}

In this section, we introduce the Jacobi--Trudi determinant
$\chi_{\lambda/\mu,a}$ of type $C_n$
and its associated paths.
Most of the information is taken from \cite{NN1}.
See~\cite{NN1} for more information.

\subsection[Jacobi-Trudi determinant of type $C_n$]{Jacobi--Trudi determinant of type $\boldsymbol{C_n}$}

A {\it partition} is a sequence of weakly decreasing
non-negative integers $\lambda=(\lambda_1, \lambda_2, \dots)$
with f\/initely many non-zero terms $\lambda_1 \ge \lambda_2 \ge\cdots \ge
\lambda_l>0$.
The {\it length} $l(\lambda)$ of $\lambda$ is
the number of the non-zero integers.
The {\it conjugate} is denoted by $\lambda'=(\lambda'_1,\lambda'_2, \dots)$.
As usual, we identify a partition $\lambda$ with a~{\it Young diagram}
$\lambda = \{ (i,j)\in \bN^2 \mid 1 \le j \le \lambda_i \}$,
and also identify a pair of partitions such that $\lambda_i \ge \mu_i$
for any $i$, with a {\it skew diagram}
$\lambda/\mu = \{ (i,j) \in \bN^2 \mid \mu_i+1 \le j \le \lambda_i\}$.

Let
\begin{equation}\label{eq:entries}
I=\{1, 2, \dots, n, \overline{n}, \dots,\overline{2},  \overline{1} \}.
\end{equation}
Let $\cZ$ be the commutative ring over $\bZ$ generated by
$z_{i,a}$'s, $i \in I$, $a \in \bC$, with the following generating relations:
\begin{gather}\label{eq:relations}
z_{i,a}z_{\overline{i},a-2n+2i-4} = z_{i-1,a}z_{\overline{i-1},a-2n+2i-4}
\quad (i=1, \dots, n), \qquad
z_{0,a}=z_{\overline{0},a}=1.
\end{gather}

Let $\bZ[[X]]$ be the formal power series ring
over $\bZ$ with variable $X$.
Let $\cA$ be the {\it non-commutative} ring
generated by $\cZ$ and $\bZ[[X]]$ with relations
\begin{equation*}
Xz_{i,a}=z_{i,a-2}X, \qquad i \in I, a \in \bC.
\end{equation*}
For any $a \in \bC$, we def\/ine $E_a(z,X)$, $H_a(z,X)\in \cA$
as
\begin{gather}
E_a(z,X) := \label{eq:E}
\left\{ \prod_{1 \le k \le n}^{\rightarrow} (1 + z_{k,a}X) \right\}
(1-z_{n,a}Xz_{\overline{n},a}X)
\left\{ \prod_{ 1 \le k \le n}^{\leftarrow}
(1 + z_{\overline{k},a}X) \right\}, \\
H_a(z,X)  := \label{eq:H}
\left\{ \prod_{1 \le k \le n}^{\rightarrow}
(1 - z_{\overline{k},a}X)^{-1} \right\}
(1 - z_{n,a}Xz_{\overline{n},a}X)^{-1}
\left\{\prod_{1 \le k \le n}^{\leftarrow} (1 - z_{k,a}X)^{-1} \right\},
\end{gather}
where $\prod\limits_{1 \le k \le n}^{\rightarrow} A_k=A_1\cdots A_n$
and $\prod\limits _{1 \le k \le n}^{\leftarrow} A_k=A_n\cdots A_1$.
Then we have
\begin{equation}\label{eq:HE}
H_a(z, X)E_a(z,-X) = E_a(z,-X)H_a(z,X)=1.
\end{equation}

For any $i \in \bZ$ and $a \in \bC$,
we def\/ine $e_{i,a}$, $h_{i,a}\in \cZ$ as
\begin{equation*}
E_a(z,X) = \sum_{i=0}^{\infty}e_{i,a}X^i, \qquad
H_a(z,X) = \sum_{i=0}^{\infty}h_{i,a}X^i,
\end{equation*}
with $e_{i,a}=h_{i,a}=0$ for $i < 0$.
Note that  $e_{i,a}=0$ if $i>2n+2$ or $i=n+1$.
It is worth mentioning that
the equality
 \begin{equation}
 e_{2n+2-i,a}=-e_{i,a-2n+2i-2}
 \end{equation}
 holds for any $i$ (cf.
\cite[equation~(2.14)]{KOSY}), though
we do not use it in the rest of the paper.
This follows from the following pseudo-antisymmetric property of $E_a(z,X)$,
\begin{equation}
E_a(z,X^{-1})\vert_{z_{i,a}\mapsto z_{\overline{i},a},
 z_{\overline{i},a}\mapsto z_{i,a}} X^{2n+2}
  = - E_{a+2n}(z,X),
\end{equation}
which is proved by successive applications of the relations (\ref{eq:relations}).

Due to the relation \eqref{eq:HE}, we have \cite[equation~(2.9)]{M}
\begin{equation}\label{eq:determinant}
\det (h_{\lambda_i-\mu_j-i+j, a+2(\lambda_i - i)})_{1 \le i,j \le l}
= \det(e_{\lambda'_i-\mu'_j-i+j, a-2(\mu'_j-j+1)})_{1 \le i,j \le l'}
\end{equation}
for any pair of partitions $(\lambda, \mu)$, where $l$ and $l'$ are
any non-negative integers such that $l \ge l(\lambda), l(\mu)$
and $l' \ge l(\lambda'), l(\mu')$.
For any skew diagram $\lambda/\mu$,
let $\chi_{\lambda/\mu,a}$ denote the determinant on the left- or
right-hand side of \eqref{eq:determinant}.
We call it the {\it Jacobi--Trudi determinant} associated with
the quantum af\/f\/ine algebra $\Uqhg$ of type $C_n$.

Let $d(\lambda/\mu):=\max \{ \lambda'_i-\mu'_i\}$ be the {\it depth}
of $\lambda/\mu$. We conjecture that, if $d(\lambda/\mu) \le n$,
the determinant $\chi_{\lambda/\mu, a}$
is the {\it $q$-character} for a certain f\/inite-dimensional representation
$V$ of the quantum af\/f\/ine algebra of type $C_n$.
If $\lambda/\mu$ is connected,
we further expect that $\chi_{\lambda/\mu, a}$ is
the $q$-character for the {\it irreducible} representation $V$
whose  highest weight monomial  is
\begin{equation}
\prod_{j=1}^{l(\lambda')}
Y_{\lambda'_j-\mu'_j,q^{a+2j-\lambda'_j-\mu'_j-1}},
\end{equation}
where $Y_{i,a}$ is the variables for the $q$-characters in \cite{FR}.

The rest of the paper is devoted to providing the combinatorial
description of the determinant $\chi_{\lambda/\mu,a}$.

\subsection[Gessel-Viennot paths]{Gessel--Viennot paths}\label{sec:GV}

Following \cite{NN1}, let us
apply the Gessel--Viennot path method (see \cite[Chapter 4.5]{S}
for good exposition)
to the determinant
$\chi_{\lambda/\mu,a}$ in \eqref{eq:determinant}
and the generating function $H_a(z, X)$ in \eqref{eq:H}.
The function $H_a(z, X)$ is the generating
function of the symmetric polynomials.
So, one can employ the $h$-labeling of \cite{S}.
However, the middle factor produces the
powers of the  pair of variables $z_{n,a}z_{\overline{n},a-2}$,
which is the source of the whole complexity \cite{NN1}.

\begin{figure}
\centerline{\includegraphics{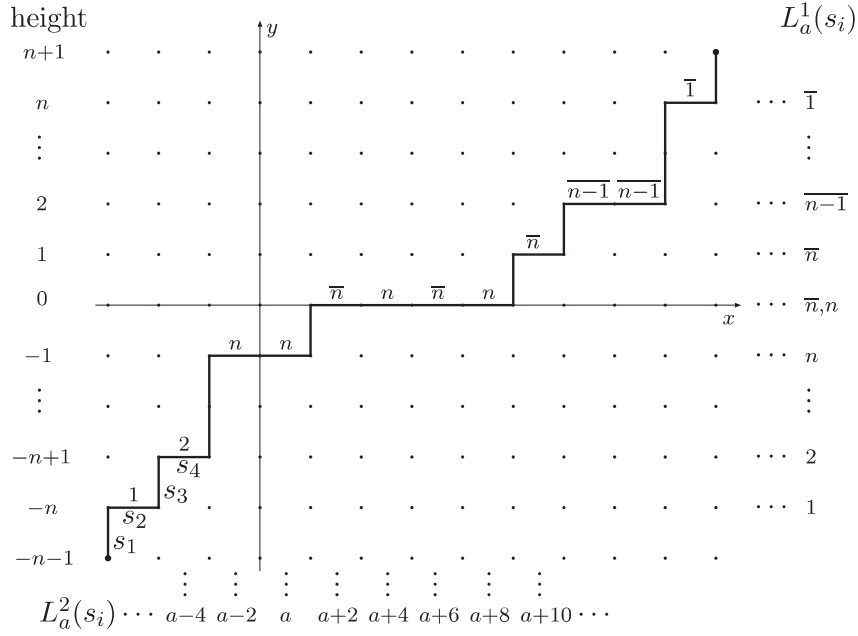}}
\caption{An example of a path of type $C_n$ and its $h$-labeling.}\label{fig:path}
\end{figure}

Consider the lattice $\bZ^2$.
An {\it E-step} (resp.\ {\it N-step})
 $s$ is a step from a point $u$ to a point $v$ in
the lattice of unit length  in east (resp.\ north) direction.
%We write these steps as $\pathNakanishi{u}{s}{v}$.
For any point $(x,y)\in \bR^2 $,
we def\/ine the {\it height} by $y$
and the {\it horizontal position} by $x$.
Due to \eqref{eq:H}, we def\/ine
a {\it path $p$} (of type $C_n$)
as a~sequence of consecutive steps
$(s_1, s_2, \dots )$ which satisf\/ies the following conditions:
\begin{enumerate}\itemsep=0pt
\item[(i)]
It starts from a point $u$ at height $-n-1$ and ends at a point $v$
at height $n+1$.
\item[(ii)]
Each step $s_i$ is an E-
or N-step.
\item[(iii)]  \label{item:D-path-three}
The number of E-steps at height 0 is even.
\item[(iv)]
The E-steps do not occur at
height $\pm(n+1)$,
\end{enumerate}
We also write such a path $p$ starting from $u$ and ending at $v$ as $\pathNakanishi{u}{p}{v}$.
See Fig.~\ref{fig:path} for an example.

\begin{rem}
We slightly change the def\/inition of paths from that of  \cite{NN1}.
In \cite{NN1}, a path starts at height $-n$ and end at height $n$,
satisfying (ii) and (iii).
Here, we extend it by adding N-steps at its both ends.
The reason of this change is  to make the formulation in
Section \ref{sec:positive-sum} as parallel as
possible to the $D_n$ case.
With our new def\/inition, a {\it boundary II-unit} in
Def\/inition \ref{def:unit}  makes sense,
and so does  a {\it II-region} in Def\/inition~\ref{def:region}.
\end{rem}

Let $\cP$ be the set of all the paths.
For any $p \in \cP$, set
\begin{gather}
E(p) := \{s \in p \mid
\text{$s$ is an  E-step} \}, \nonumber\\
E_0(p) := \{s \in p \mid
\text{$s$ is an E-step at height 0} \} \subset E(p).\label{eq:label-steps}
\end{gather}
If $E_0(p)=\{s_j, s_{j+1}, \dots, s_{j+2k-1} \}$, then
let
\begin{equation*}
E_0^1(p):=\{s_{j+1}, s_{j+3}, \dots, s_{j+2k-1} \}
\subset E_0(p).
\end{equation*}
Fix $a \in \bC$. The {\it $h$-labeling} (of type $C_n$) associated
with $a$ for a path $p \in \cP$ is the pair of maps
$L_a = (L^1_a, L^2_a)$ on $E(p)$ def\/ined as follows:
Suppose that a step $s\in E(p)$ starts at a point $w=(x,y)$.
Then, we set{\samepage
\begin{gather}\label{eq:e-label}
L_a^1(s)  =
\begin{cases}
n+1+y, & \text{if $y < 0$}, \\
\overline{n}, & \text{if $y=0$ and $s \not\in E_0^1(p)$}, \\
n, & \text{if $y=0$ and $s \in E_0^1(p)$}, \\
\overline{n+1-y}, & \text{if $y > 0$}, \\
\end{cases}
\\
L_a^2(s) = a + 2x.\nonumber
\end{gather}
See Fig.~\ref{fig:path}.}

Now we def\/ine the {\it weight} of $p \in \cP$ as
\begin{equation*}
z_a^p := \prod_{s \in E(p)} z_{L_a^1(s), L_a^2(s)} \in \cZ.
\end{equation*}
By the def\/inition of $E_a(z,X)$ in \eqref{eq:E},
for any $k\in \bZ$, we have
\begin{equation}\label{eq:e}
h_{r, a+2k+2r-2}(z) = \sum_pz_a^p,
\end{equation}
where the sum runs over all $p \in \cP$
such that $(k, -n-1) \overset{p}{\to} (k+r, n+1)$.

For any $l$-tuples of distinct points
$\bu = (u_1, \dots, u_l)$ of height $-n-1$
and $\bv = (v_1, \dots, v_l)$ of height $n+1$,
and any permutation $\sigma \in \fS_l$, let
\begin{gather*}
\fP(\sigma; \bu, \bv)= \{ \bp = (p_1, \dots, p_l) \mid
p_i\in \cP,\ \pathNakanishi{u_i}{p_i}{v_{\sigma(i)}} \text{ for } i=1, \dots, l \},
\end{gather*}
and set
\begin{gather*}
\fP(\bu, \bv) = \bigsqcup_{\sigma \in \fS_l}\fP(\sigma; \bu, \bv).
\end{gather*}
We def\/ine the {\it weight} $z_a^{\bp}$ and the {\it sign} $(-1)^{\bp}$ of
$\bp \in \fP(\bu, \bv)$
as
\begin{equation}\label{eq:weight}
z_a^{\bp} := \prod_{i=1}^lz_a^{p_i}, \qquad
(-1)^{\bp}:= \sgn \, \sigma \quad \text{if $\bp \in \fP(\sigma; \bu, \bv)$}.
\end{equation}

For any skew diagram $\lambda/\mu$, set $l=l(\lambda)$, and
\begin{align*}
\bu_{\mu}=(u_1, \dots, u_{l}), \qquad
& u_i=(\mu_i + 1 -i, -n -1),  \\
\bv_{\lambda}=(v_1, \dots, v_{l}), \qquad
& v_i=(\lambda_i + 1 -i, n +1).
\end{align*}
Then, due to \eqref{eq:e}, the determinant \eqref{eq:determinant}
can be written as
\begin{equation}\label{eq:GV-path-description}
\chi_{\lambda/\mu,a}
= \sum_{\bp \in \fP(\bu_{\mu}, \bv_{\lambda})} (-1)^{\bp}z_a^{\bp}.
\end{equation}

\begin{defn}
We say that an intersecting pair $(p_i,p_j)$ of paths is
{\it specially intersecting} if it satisf\/ies the following conditions:
\begin{enumerate}\itemsep=0pt
\item The intersection of $p_i$ and $p_j$ occurs only at height $0$.
\item $p_i(0)-p_j(0)$ is odd, where
$p_i(0)$ is the horizontal position of
the leftmost point on $p_i$ at height $0$.
\end{enumerate}
Otherwise, we say that
an intersecting pair $(p_i,p_j)$ is {\it ordinarily intersecting}.
\end{defn}

We can def\/ine a weight-preserving, sign-reversing involution
on the set of all the tuples $\bp \in \fP(\bu_{\mu}, \bv_{\lambda})$
which have some ordinarily intersecting pairs $(p_i, p_j)$.
Therefore, we have
\begin{prop}[{\cite[Proposition 5.3]{NN1}}]\label{prop:first-involution}
For any skew diagram $\lambda/\mu$,
\begin{equation}\label{eq:first-sum}
\chi_{\lambda/\mu, a}= \sum_{\bp \in P_1(\lambda/\mu)}(-1)^{\bp}z_a^{\bp},
\end{equation}
where $P_1(\lambda/\mu)$ is the set of all
$\bp \in \fP(\bu_{\mu}, \bv_{\lambda})$ which do not have any ordinarily
intersecting pair of paths.
\end{prop}

\section{Paths description }\label{sec:positive-sum}
In this section, we turn the expression \eqref{eq:first-sum}
 into an expression  (\ref{eq:path-description})
 by a positive sum over a set of tuples of paths
(`paths description'), which will be naturally translated to
 the tableaux description Section \ref{sec:tableau-description}.

 \subsection{Outline of formulation}
 As we mentioned in the introduction,
the formulation here is quite parallel to the $D_n$ case \cite{NN2}.
Before starting, let us brief\/ly explain the idea behind the formulation.

The expression \eqref{eq:first-sum} still involves negative terms.
Since $\chi_{\lambda/\mu,a}$ is supposed to be the
$q$-character of a certain representation if $d(\lambda/\mu)\leq n$,
these negative terms should be canceled by
the positive terms of the same weight if $d(\lambda/\mu)\leq n$.
(It turned out that the condition to have a positive
 sum can be slightly extended to
(\ref{eq:positivity}).)
This motivates us to construct a new type of
weight-preserving, sign-reversing
involution besides the standard one (the {\it first involution} in~\cite{NN2})
 used for
 Proposition \ref{prop:first-involution}.
 A key  is  the relations (\ref{eq:relations})
among our variables, which allow to deform a tuple of paths preserving its weight.
 This deformation is called the {\it expansion} and {\it folding} in \cite{NN2}.
 Using them, one can construct the desired involution $\iota_2$  (the {\it second involution}
 in \cite{NN2}), and obtain an expression (\ref{eq:positive-sum}) by a positive sum.

 However, it turns out that (\ref{eq:positive-sum}) is not our f\/inal answer yet,
 because a tuple of paths for~(\ref{eq:positive-sum}) may be transposed (i.e.,
 there is a pair of paths whose orders of the initial points and f\/inal points are
 transposed); in that case it cannot be  translated it into a tableau
 on $\lambda/\mu$ \`a la Gessel--Viennot~\cite{GV, S}.
 To resolve this problem, we construct a weight-preserving map $\phi$
 (the {\it folding map} in~\cite{NN2}), which transforms a tuple of paths
 of (\ref{eq:positive-sum}) into non-transposed one.
 Then, we f\/inally obtain the desired expression (\ref{eq:path-description})
 which is the counterpart of   the tableaux description in Section \ref{sec:tableau-description}.

 Below we are going to def\/ine basic objects
 (lower/upper path, I/II-unit, even/odd I/II-region, etc.),
all of which are introduced to construct the above mentioned
maps $\iota_2$ and $\phi$.
 The dif\/ferences between $C_n$ and $D_n$,
 which originate from the ones of the generating functions
 and the relations of the variables,
  are absorbed into
 the def\/initions of these objects.
 With these objects,
 the construction of $\iota_2$ and $\phi$ are
done in a unif\/ied way for  both $C_n$ and $D_n$.
For the reader who is only interested in the f\/inal result (\ref{eq:path-description}),
only the def\/inition of an {\it odd II-region} is crucial.
See Fig.~\ref{fig:I-region} to have  idea of a II-region.

\subsection[$\I$- and $\II$-regions]{$\boldsymbol{\I}$- and $\boldsymbol{\II}$-regions}
\begin{figure}[t]
\centerline{\includegraphics{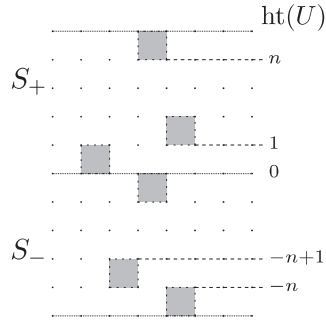}}
\caption{Examples of units.}\label{fig:units}
\end{figure}

\begin{figure}[t]
\centerline{\includegraphics{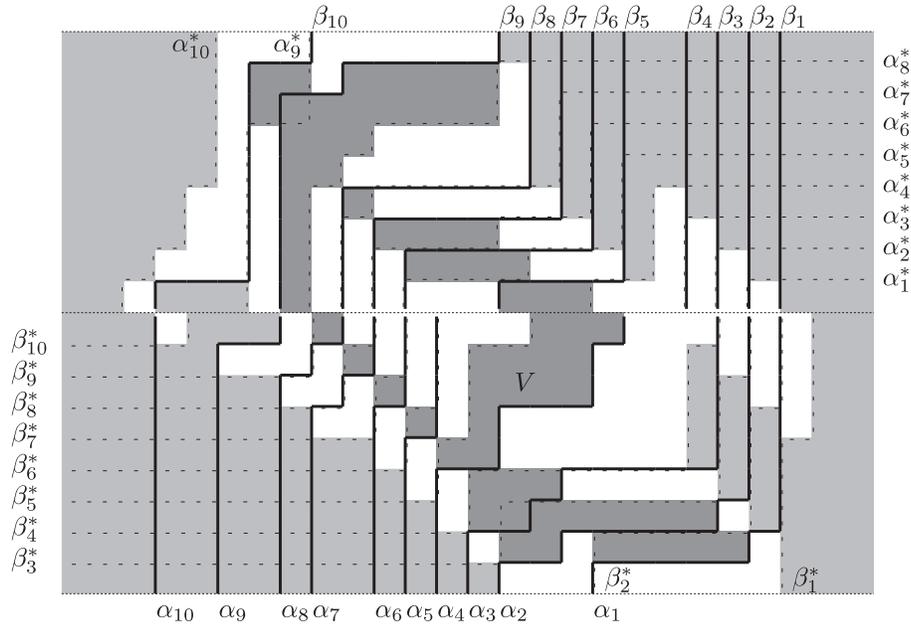}}
\caption{An example of $(\alpha; \beta)$ and its
dual $(\alpha^*; \beta^*)$.
The $\II$-units of $(\alpha; \beta)$ are shaded, and,
especially, the darkly shaded region
represents a $\II$-region $V$.}
\label{fig:I-region}
\end{figure}

Here, we introduce some
notions which are necessary to give the expression
by a positive sum.
See Figs.~\ref{fig:units} and \ref{fig:I-region} for examples.

Let
\begin{equation}
S_+  := \bR \times [0,n+1],\qquad
S_-  := \bR \times [-n-1,0].
\end{equation}

\begin{defn}
A {\it lower path} $\alpha$ (of type $C_n$)
is a sequence of consecutive steps in $S_-$ which starts at a
point of height $-n-1$ and
ends at a point of height $0$ such that each step is
an E- or N-step, and an E-step does not occur at
height $0$ and $-n-1$.
Similarly, an {\it upper path} $\beta$ (of type $C_n$)
is a sequence of consecutive steps in $S_+$ which
starts at a point of height $0$ and
ends at a point of height $n+1$ such that each step is
an E- or N-step, and an E-step does not occur at
height $0$ and $n+1$.
\end{defn}

For a lower path $\alpha$ and an upper path $\beta$,
let $\alpha(r)$ and $\beta(r)$ be the horizontal positions
of the leftmost points of $\alpha$ and
the rightmost points $\beta$ at height $r$, respectively.

Let
\begin{equation}
(\alpha; \beta):=(\alpha_1, \dots, \alpha_l; \beta_1, \dots, \beta_l)
\end{equation}
be a pair of an $l$-tuple $\alpha$
of lower paths and an $l$-tuple $\beta$ of upper paths.
We say that $(\alpha; \beta)$ is {\it nonintersecting}
if $(\alpha_i, \alpha_j)$ is not intersecting, and so is
$(\beta_i, \beta_j)$ for any $i \ne j$.

{}From now on,
let $\lambda/\mu$ be a skew diagram, and we set $l=l(\lambda)$.
Let
\begin{gather*}
\cH(\lambda/\mu):=
\Bigg\{ (\alpha; \beta)=(\alpha_1, \dots, \alpha_l; \beta_1, \dots, \beta_l)
\Bigg|
\begin{array}{l}
(\alpha; \beta) \text{ is nonintersecting}, \\
\alpha_i(-n-1)=\mu_i+1-i, \\
\beta_i(n+1)=\lambda_i+1-i
\end{array}
\Bigg\}.
\end{gather*}
For any  $(\alpha;\beta)\in \cH(\lambda/\mu)$,
we def\/ine $\alpha_i^*(r)$, $\beta_i^*(-r)$,
($i=1$, \dots, $l$; $r=0$, \dots, $n$) as
\begin{gather}\label{eq:alpha-dual}
\alpha_i^*(r)=
\begin{cases}
\alpha_{i-r}(-r)-1-r& i-r\geq 1,\\
+\infty & i-r \leq 0,
\end{cases}
\\
\label{eq:beta-dual}
\beta_i^*(-r)=
\begin{cases}
\beta_{i+r}(r)+1+r& i+r\leq l,\\
-\infty & i+r \geq l+1.
\end{cases}
\end{gather}
Since $(\alpha;\beta)$ is nonintersecting,
we have $\alpha_i^*(r)\leq \alpha_i^*(r+1)$
and $\beta_i^*(-r)\leq \beta_i^*(-r+1)$.
Therefore, one can naturally identif\/ies the
data $\{ \alpha_i^*(r)\}_{r=0}^n$ with
an upper path $\alpha_i^*$ (an `inf\/inite' upper path
if it contains $+\infty$).
The lower path $\beta_i^*$ is obtained similarly.
It is easy to check that $(\alpha^*;\beta^*)$ is nonintersecting.
We call $(\alpha^*;\beta^*)$ the {\it dual} of
$(\alpha,\beta)$.
As we mention in the introduction,
the def\/inition of $(\alpha^*;\beta^*)$ is more delicate
than the one for $D_n$,
where $\alpha_i^*(r):=\alpha_i(-r)-1$ and $\beta^*_i(-r):=
\beta_i(r)+1$.
It is so def\/ined that the involution $\iota_2$
in Proposition \ref{prop:second-involution} is
weight-preserving under the relation (\ref{eq:relations}).

For any skew diagram $\lambda/\mu$,
we call the following condition the {\it positivity condition}:
\begin{equation}\label{eq:positivity}
d(\lambda/\mu)\leq n+1,
\end{equation}
or equivalently, $\lambda_{i+n+1}\leq \mu_i$ for any $i$.
We call this the `positivity condition',
because \eqref{eq:positivity} guarantees that
$\chi_{\lambda/\mu,a}$ is a positive (nonnegative, strictly to say) sum
(see Theorem \ref{thm:positive-sum}).
By the def\/inition, we have
\begin{lem}[cf.\ {\cite[Lemma 4.2]{NN2}}]\label{lem:intersecting}
Let $\lambda/\mu$ be a skew diagram satisfying
the positivity condition \eqref{eq:positivity}, and
let $(\alpha; \beta) \in \cH(\lambda/\mu)$.
Then,
\begin{equation}
\beta_{i+1}(n+1) \le \alpha^*_i(n+1), \qquad
\beta^*_{i+1}(-n-1) \le \alpha_i(-n-1).
\end{equation}
\end{lem}

A {\it unit} $U\subset S_{\pm}$ is
a unit square with its vertices on the lattice.
The {\it height} $\Ht(U)$ of $U\subset S_+$ (resp. $U \subset S_-$)
is given by the height of the
lower edge (resp. the upper edge) of $U$.

\begin{defn}\label{def:unit}
Let $(\alpha; \beta) \in \cH(\lambda/\mu)$.
For any unit $U \subset S_{\pm}$, let $\pm r=\Ht(U)$ and
let $a$ and $a'=a+1$ be the horizontal positions of
the left and the right edges of $U$.
Set $\beta_{l+1}(r)=\beta^*_{l+1}(-r)=-\infty$
and $\alpha_0(-r)=\alpha^*_0(r)=+\infty$.
Then
\begin{enumerate}\itemsep=0pt
\item $U$ is called a {\it $\I$-unit} of $(\alpha; \beta)$ if
there exists some $i$ ($0 \le i \le l$) such that
\begin{gather}
 \alpha^*_i(r) \le a < a' \le \beta_{i+1}(r), \qquad
\text{if $U \subset S_+$}, \nonumber\\
 \alpha_i(-r) \le a < a' \le \beta^*_{i+1}(-r), \qquad
\text{if $U \subset S_-$}.\label{eq:I-unit}
\end{gather}
\item $U$ is called a {\it $\II$-unit} of $(\alpha; \beta)$ if
there exists some $i$ ($0 \le i \le l$) such that
\begin{gather}
 \beta_{i+1}(r) \le a < a' \le \alpha_i^*(r), \qquad
\text{if $U \subset S_+$}, \nonumber\\
 \beta^*_{i+1}(-r) \le a < a' \le \alpha_i(-r), \qquad
\text{if $U \subset S_-$}.\label{eq:II-unit}
\end{gather}
\end{enumerate}
Furthermore, a $\II$-unit $U$ of $(\alpha; \beta)$
is called a {\it boundary} $\II$-unit if
one of the following holds:
\begin{itemize}\itemsep=0pt
\item[(i)] $U\subset S_+$, and  \eqref{eq:II-unit}
holds for $i\leq r $, $i=l$, or $r=n$.
\item[(ii)] $U\subset S_-$, and  \eqref{eq:II-unit}
holds for $i=0 $, $i\geq l-r$, or $r=n$.
\end{itemize}
\end{defn}

For a unit $U\subset S_+$ with
 vertices $(x,y)$, $(x+1,y)$, $(x,y+1)$, $(x+1,y+1)$,
the {\it dual\/} $U^*\subset S_-$ is a unit
with vertices
$(x+1+y,-y)$, $(x+2+y,-y)$, $(x+1+y,-y-1)$, $(x+2+y,-y-1)$.
Conversely, we def\/ine $(U^*)^*=U$.
Let $U$ and $U'$ be units.
If the upper-left or the lower-right vertex of $U$
is also a vertex of $U'$, then
we say that $U$ and $U'$ are {\it adjacent}
and write $U\diamond U'$.

So far, the def\/initions are specif\/ic for type $C_n$.
{}From now till the end of Section \ref{sec:positive-sum},
all the statements are
literally the same as type $D_n$ \cite{NN2}.

Fix $(\alpha; \beta)\in \cH(\lambda/\mu)$.
Let $\cU_{\I}$ be the set of all $\I$-units of $(\alpha; \beta)$,
and let
$\tilde{\cU}_{\I}:=\bigcup_{U\in \cU_{\I}}U$,
where the union is taken for $U$ as a subset of $S_+ \sqcup S_-$.
Let $\sim$ be the equivalence relation in $\cU_{\I}$ generated by
the relation $\diamond$, and $[U]$ be its equivalence class of $U \in \cU_{\I}$.
We call $\bigcup_{U'\in [U]}U'$
a {\it connected component} of $\tilde{\cU}_{\I}$.
For $\II$-units,
$\cU_{\II}$, $\tilde{\cU}_{\II}$ and its connected component
are def\/ined similarly.

Now we introduce the main concept in the section.
\begin{defn}
\label{def:region}
Let $\lambda/\mu$ be a skew diagram satisfying the positivity condition
\eqref{eq:positivity},
and let $(\alpha; \beta) \in \cH(\lambda/\mu)$.
\begin{enumerate}\itemsep=0pt
\item
A connected component $V$ of $\tilde{\cU}_{\I}$
is called a {\it $\I$-region} of $(\alpha; \beta)$
if it contains at least one $\I$-unit of height $0$.
\item
A connected component $V$ of $\tilde{\cU}_{\II}$
is called a {\it $\II$-region} of $(\alpha; \beta)$
if it satisf\/ies the following conditions:
\begin{enumerate}[(i)]\itemsep=0pt
\item $V$ contains at least one $\II$-unit of height $0$.
\item $V$ does not contain any boundary $\II$-unit.
\end{enumerate}
\end{enumerate}
\end{defn}

\begin{prop}[cf.\ {\cite[Proposition 4.6]{NN2}}]\label{prop:duality}
If $V$ is a $\I$- or  $\II$-region, then
$V^*=V$,
where
for a union of units $V=\bigcup U_i$,
we define $V^*=\bigcup U_i^*$.
\end{prop}

\subsection{Second involution}
\label{subsec:positive-sum}

Following the $D_n$ case \cite{NN2},
let us  derive an  expression
of $\chi_{\lambda/\mu,a}$ by a positive sum
from \eqref{eq:first-sum}.

Since the proofs for $D_n$ are applicable almost word for word
to all the statements below, we omit them.
For interested readers, we provide some technical information
in Appendix \ref{sec:appendix}.

{}From now on, we assume that
$\lambda/\mu$ satisf\/ies the positivity condition \eqref{eq:positivity}.
For each $\bp \in P_1(\lambda/\mu)$, one can uniquely
associate $(\alpha; \beta) \in \cH(\lambda/\mu)$
by removing all the E-steps from $\bp$.
We denote
by $\pi(p)$ the element $(\alpha;\beta)$ corresponding to the path $p$.
A $\I$- or $\II$-region of $(\alpha; \beta)=\pi(\bp)$
is also called a {\it $\I$- or $\II$-region of $\bp$}.
If $h:=\alpha_i(0)-\beta_{i+1}(0)$ is
a non-positive number (resp.\ a~positive number),
then we call a pair $(\alpha_i, \beta_{i+1})$ an {\it overlap}
(resp.~a {\it hole}).
Furthermore, if $h$ is an even number (resp.~an odd number),
then we say that $(\alpha_i, \beta_{i+1})$ is {\it even} (resp.~{\it odd}).
For any $\I$-region $V$ (resp.~$\II$-region $V$) of
 $p\in P_1(\lambda/\mu)$ with $(\alpha;\beta)=\pi(p)$, we set
\begin{equation}\label{eq:m}
n(V):=\# \Big\{ i \, \Big| \text{\begin{tabular}{l}
$(\alpha_i, \beta_{i+1})$
is an even overlap (resp.~an even hole)\\
which intersects with $V$ at height $0$
\end{tabular}}\Big\}.
\end{equation}
For example, $n(V)=2$ for $V$ in Fig.~\ref{fig:I-region}.

\begin{defn}
We say that a $\I$- or $\II$-region $V$ is
{\it even} (resp.~{\it odd}) if $n(V)$ is even (resp.~odd).
\end{defn}

Let
$P_{\odd}(\lambda/\mu)$ be the set of all
$\bp \in P_1(\lambda/\mu)$ which have at least one
odd $\I$- or $\II$-region of $\bp$.

\begin{prop}[cf.\ {\cite[Proposition 4.12]{NN2}}]
\label{prop:second-involution}
There exists
a weight-preserving, sign-rever\-sing involution
$\iota_2: P_{\odd}(\lambda/\mu) \to P_{\odd}(\lambda/\mu)$.
\end{prop}

It follows from Proposition \ref{prop:second-involution} that
the contributions of $P_{\odd}(\lambda/\mu)$
to the sum \eqref{eq:first-sum} cancel each other.
Let
$P_2(\lambda/\mu):=
P_1(\lambda/\mu) \backslash P_{\odd}(\lambda/\mu)$, i.e.,
the set of all
$\bp \in P_1(\lambda/\mu)$ which satisfy the following conditions:
\begin{enumerate}[(i)]
\item \label{item:Ptwo-no-ord}
$\bp$ does not have any ordinarily intersecting pair $(p_i, p_j)$.
\item \label{item:Ptwo-no-odd}
$\bp$ does not have any odd $\I$- or $\II$-region.
\end{enumerate}
For any $\bp \in P_2(\lambda/\mu)$,
$(-1)^{\bp}=1$ holds.
Thus, the sum \eqref{eq:first-sum} reduces to a positive sum,
and we have an  expression by a positive sum,
\begin{thm}[cf.\  {\cite[Theorem 4.13]{NN2}}]
\label{thm:positive-sum}
For any skew diagram $\lambda/\mu$
satisfying the positivity condition \eqref{eq:positivity},
we have
\begin{equation}\label{eq:positive-sum}
\chi_{\lambda/\mu, a}= \sum_{\bp \in P_2(\lambda/\mu)}z_a^{\bp}.
\end{equation}
\end{thm}

\subsection{Paths description}
Since a tuple of paths $\bp  \in \fP(\sigma; \bu_{\mu}, \bv_{\lambda})$ is naturally translated into
a tableau of shape $\lambda/\mu$
if and only if $\sigma = \id $,
we introduce another set of paths as follows.
Let $P(\lambda/\mu)$ be the set of all
$\bp \in \fP(\id; \bu_{\mu}, \bv_{\lambda})$
such that
\begin{enumerate}[(i)]\itemsep=0pt
\item $p$ does not have any ordinarily intersecting
{\it adjacent} pair $(p_i, p_{i+1})$.
\item $p$ does not have any odd $\II$-region.
\end{enumerate}
Here, an odd $\II$-region of $p\in P(\lambda/\mu)$ is
def\/ined in the same way as that of $p\in P_1(\lambda/\mu)$.
The following fact is not so trivial.
\begin{prop}
[cf.\ {\cite[Proposition 5.1]{NN2}}]
\label{prop:folding-map}
There exists a weight-preserving bijection
\begin{equation*}
\phi: P_2(\lambda/\mu) \to P(\lambda/\mu).
\end{equation*}
\end{prop}
The map $\phi$ is called the {\it folding map} in \cite{NN2}.
{}From Theorem \ref{thm:positive-sum} and Proposition \ref{prop:folding-map},
we immediately have
\begin{thm}
[Paths description, cf.\ {\cite[Theorem 5.2]{NN2}}]
\label{thm:tilde-positive-sum}
For any skew diagram $\lambda/\mu$ satisfying
the positivity condition \eqref{eq:positivity},
we have
\begin{equation}\label{eq:path-description}
\chi_{\lambda/\mu, a}= \sum_{\bp \in P(\lambda/\mu)} z_a^{\bp}.
\end{equation}
\end{thm}

\section{Tableaux description}
\label{sec:tableau-description}

\subsection{Tableaux description}\label{sec:extra}
Def\/ine a total order in $I$ in \eqref{eq:entries} by
\begin{equation*}
1 \prec 2 \prec \dots \prec
n \prec \overline{n}
\prec \dots \prec \overline{2} \prec\overline{1}.
\end{equation*}

A {\it tableau} $T$ of shape $\lambda/\mu$ is the skew diagram $\lambda/\mu$
with each box f\/illed by one entry of $I$.
For a tableau $T$ and $a \in \bC$, we def\/ine
the {\it weight} of $T$ as
\begin{equation*}
z_a^T=\prod_{(i,j) \in \lambda/\mu}z_{T(i,j),a+2(j-i)},
\end{equation*}
where $T(i,j)$ is the entry of $T$ at $(i,j)$.

\begin{defn}
A tableau $T$ (of shape $\lambda/\mu$) is
called an {\it HV-tableau} if it satisf\/ies the following conditions:
\begin{itemize}\itemsep=0pt
\item[$(\bH )$]  horizontal rule.
Each $(i,j)\in \lambda/\mu$ satisf\/ies both of the following conditions:
\begin{itemize}\itemsep=0pt
\item[(i)] $T(i,j)\preceq T(i,j+1)$ or $(T(i,j), T(i,j+1))=(\overline{n}, n)$.
\item[(ii)] $(T(i,j-1), T(i,j), T(i,j+1))\ne
    (\overline{n}, \overline{n}, n), (\overline{n}, n, n)$.
\end{itemize}
\item[$(\bV )$] vertical rule.
Each $(i,j)\in \lambda/\mu$ satisf\/ies
one of the following conditions:
\begin{itemize}\itemsep=0pt
\item[(i)] $T(i,j)\prec T(i+1,j)$.
\item[(ii)] $T(i,j)=T(i+1,j)=n$,
$(i+1,j-1)\in \lambda/\mu$, $T(i+1,j-1)=\overline{n}$.
\item[(iii)] $T(i,j)=T(i+1,j)=\overline{n}$,
$(i,j+1)\in \lambda/\mu$, $T(i,j+1)=n$.
\end{itemize}
\end{itemize}
\end{defn}

The rule $(\bH)$ appears in \cite{KS} for a skew diagram
of one row.
We write the set of all HV-tableaux of shape $\lambda/\mu$
by $\HVTab(\lambda/\mu)$.

Let $\HVP(\lambda/\mu)$ be the set of all
$\bp\in \fP(\id; \bu_{\mu}, \bv_{\lambda})$
which do not have any ordinarily intersecting
adjacent pair $(p_i, p_{i+1})$.
With any $\bp \in \HVP(\lambda/\mu)$, we associate a tableau $T$
of shape $\lambda/\mu$ as follows:
For any $j=1, \dots, l$,
let $E(p_j)= \{ s_{i_1}, s_{i_2}, \dots, s_{i_m} \}$
$(i_1<i_2< \dots<i_m)$
be the set def\/ined as in \eqref{eq:label-steps}, and set
\begin{equation*}
T(j,\mu_j+k)=L_a^1(s_{i_k}), \qquad k=1, \dots, m,
\end{equation*}
where $L_a^1$ is the f\/irst component of the $h$-labeling \eqref{eq:e-label}.
It is easy to see that $T$ satisf\/ies the horizontal rule $(\bH)$
because of the def\/inition of the $h$-labeling,
and satisf\/ies the vertical rule~$(\bV)$ because $\bp$
does not have any ordinarily intersecting adjacent pair.
Therefore, if we set $\cTh:\bp \mapsto T$, we have
\begin{prop}[cf.\ {\cite[Proposition 5.5]{NN2}}]
\label{prop:path-tableau}
The map
\begin{equation*}
\cTh: \HVP(\lambda/\mu) \to \HVTab(\lambda/\mu)
\end{equation*}
is a weight-preserving bijection.
\end{prop}

Note that $P(\lambda/\mu)\subset  \HVP(\lambda/\mu)$.
Let $\Tab(\lambda/\mu):=\cTh(P(\lambda/\mu))$.
In other words, $\Tab(\lambda/\mu)$ is the set of
all the tableaux $T$ which satisfy
$(\bH)$, $(\bV)$, and the following {\it extra rule}:
\vspace*{5pt}

\noindent
\begin{tabular}{ll}
$(\bE)$ & The corresponding $\bp = \cTh^{-1}(T)$ does not have
any odd $\II$-region.
\vspace{5pt}
\end{tabular}
\par\noindent
By Theorem \ref{thm:tilde-positive-sum} and Proposition \ref{prop:path-tableau},
we obtain a tableaux description of $\chi_{\lambda/\mu,a}$,
which is the main result of the paper.
\begin{thm}[Tableaux description, cf.\ {\cite[Theorem 5.6]{NN2}} ]
\label{thm:tableau-description}
For any skew diagram $\lambda/\mu$ satisfying
the positivity condition \eqref{eq:positivity},
we have
\begin{equation*}
\chi_{\lambda/\mu, a} = \sum_{T \in \Tab(\lambda/\mu)} z_a^{T}.
\end{equation*}
\end{thm}

\begin{figure}[t]
\centerline{\includegraphics{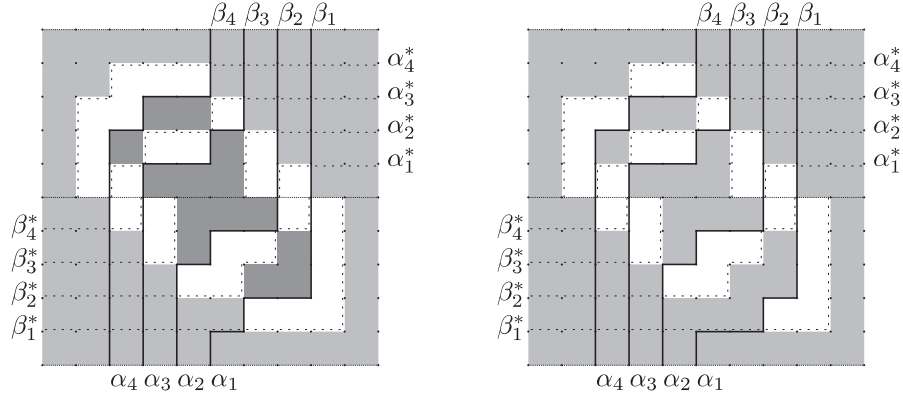}}
\caption{The tuples of paths corresponding to
$T$ and $T'$ in \eqref{eq:two-tableaux}.
The $\II$-units of $(\alpha; \beta)$ are shaded, and,
especially, the darkly shaded region
represents a $\II$-region.}\label{fig:extra-rule-example}
\end{figure}

\begin{exmp}
\label{exmp:two-tableaux}
Let $n=4$.
Consider the following two HV-tableaux which dif\/fer in only one letter:
\begin{equation}
\label{eq:two-tableaux}
T=
\raisebox{-26pt}
{
{\setlength{\unitlength}{0.2mm}
\begin{picture}(25,100)
\multiput(0,0)(25,0){4}{\line(0,1){100}}
\multiput(0,0)(0,25){5}{\line(1,0){75}}
%\multiput(0,100)(0,-25){3}{\line(1,0){25}}
\put(7,80){$1$}
\put(7,55){$3$}
\put(7,30){$\overline{4}$}
\put(7,5){$\overline{3}$}
\put(32,80){$2$}
\put(32,55){$4$}
\put(32,30){$\overline{4}$}
\put(32,5){$\overline{2}$}
\put(57,80){$2$}
\put(57,55){$4$}
\put(57,30){$\overline{3}$}
\put(57,5){$\overline{2}$}
\end{picture}
}
}
\hskip25pt
,
\quad\quad
T'=
\raisebox{-26pt}
{
{\setlength{\unitlength}{0.2mm}
\begin{picture}(25,100)
\multiput(0,0)(25,0){4}{\line(0,1){100}}
\multiput(0,0)(0,25){5}{\line(1,0){75}}
%\multiput(0,100)(0,-25){3}{\line(1,0){25}}
\put(7,80){$1$}
\put(7,55){$3$}
\put(7,30){$\overline{4}$}
\put(7,5){$\overline{3}$}
\put(32,80){$1$}
\put(32,55){$4$}
\put(32,30){$\overline{4}$}
\put(32,5){$\overline{2}$}
\put(57,80){$2$}
\put(57,55){$4$}
\put(57,30){$\overline{3}$}
\put(57,5){$\overline{2}$}
\end{picture}
}
}
\hskip25pt
.
\end{equation}
By Fig.~\ref{fig:extra-rule-example},
we see that $T\not\in\Tab(\lambda/\mu)$ because
$\cTh^{-1}(T)$ has an odd $\II$-region,
while
 $T'\in\Tab(\lambda/\mu)$ because
$\cTh^{-1}(T')$ does not so.
\end{exmp}

\subsection[Transformation into $e$-picture]{Transformation into $\boldsymbol{e}$-picture}

The paths we have used so far
are so-called `$h$-paths' in the Gessel--Viennot method.
To translate  the extra rule
$(\bE)$ into tableau language in a more explicit way,
it is convenient to transform the def\/inition of
$\II$-regions by  $h$-paths into the one by `$e$-paths'.
This procedure is not necessary for $D_n$,
where $e$-paths are employed from beginning \cite{NN2}.

\begin{defn}
A {\it lower $e$-path} $\gamma$ (of type $C_n$)
is a sequence of consecutive steps in $S_-$ which starts at a
point of height $-n-1$ and
ends at a point of height $0$ such that each step is
a W(est)- or N-step, and a W-step
occurs at most
once at each height and  does not occur at
height $0$, $-n-1$.
Similarly, an {\it upper $e$-path} $\delta$ (of type $C_n$)
is a sequence of consecutive steps in $S_+$ which
starts at a point of height $0$ and
ends at a point of height $n+1$ such that each step is
a W- or N-step, and a W-step
occurs at most
once at each height
and does not occur at
height $0$, $n+1$.
\end{defn}

For a lower $e$-path $\gamma$ and an upper $e$-path $\delta$,
let $\gamma(r)$ (resp.\ $\delta(r)$) be the horizontal position
of the rightmost point of $\gamma$ (resp. the leftmost point of
$\delta$) at height $r$.

\begin{figure}[t]
\centerline{\includegraphics{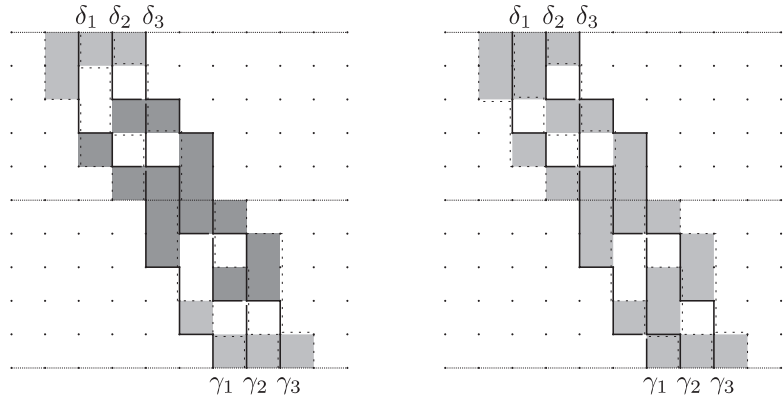}}
\caption{$(\gamma;\delta)$ corresponding to $T$ and $T'$ in \eqref{eq:two-tableaux}.
The $\II'$-units of $(\gamma; \delta)$ are shaded, and,
especially, the darkly shaded region represents
a $\II'$-region.}
\label{fig:gamma}
\end{figure}

Fix a given $p\in \HVP(\lambda/\mu)$,
and let $T=\cTh(p)$, and $l'=\lambda_1$.
With $p$ we associate
a pair of an $l'$-tuple $\gamma$
of lower $e$-paths and an $l'$-tuple $\delta$ of upper $e$-paths,
\begin{equation}
(\gamma; \delta):=(\gamma_1, \dots, \gamma_{l'};
 \delta_1, \dots, \delta_{l'}),
\end{equation}
where,
for each $i$,
we set $\gamma_i(-n-1)=i-\mu'_i$,
$\delta_i(n+1)=i-\lambda'_i$, and
we identify the $E$-steps of~$p$ at nonzero heights corresponding to
the letters in the $i$th column of $T$ with the $W$-steps
of~$\gamma_i$ and~$\delta_i$ (ignoring the direction of W- and E-steps).
We write $(\gamma; \delta)=\pi'(\bp)$.
The {\it dual} $\gamma^*_i$ of $\gamma_i$ is the upper $e$-path
with
$\gamma_i^*(r)=\gamma_i(-r)-1-r$ ($r=0,\dots,n$).
Similarly,
The {\it dual} $\delta^*_i$ of $\delta_i$ is the lower $e$-path
with $\delta_i^*(-r)=\delta_i(r)+1+r$ ($r=0,\dots,n$).
See Fig.~\ref{fig:gamma} for examples of
notions in this subsection.

\begin{defn}
\label{def:II'-unit}
Let $(\gamma; \delta)$ be as above.
For any unit $U \subset S_{\pm}$, let $\pm r=\Ht(U)$ and
let $a$ and $a'=a+1$ be the horizontal positions of
the left and the right edges of $U$.
Then,
$U$ is called a {\it $\II'$-unit} of $(\gamma; \delta)$ if
there exists some $i$ ($1 \le i \le l'$) such that
\begin{gather}
 \gamma_{i}^*(r) \le a < a' \le \delta_i(r), \qquad
\text{if $U \subset S_+$}, \nonumber\\
 \gamma_{i}(-r) \le a < a' \le \delta_i^*(-r), \qquad
\text{if $U \subset S_-$}.\label{eq:e-II-unit}
\end{gather}
Furthermore, a $\II'$-unit $U$ of $(\gamma; \delta)$
is called a {\it boundary} $\II'$-unit if
\eqref{eq:e-II-unit} holds for $r= n$.
\end{defn}

A connected component
of the union $\tilde{\cU}_{\II'}$
of all $\II'$-units of $(\gamma; \delta)$
is def\/ined similarly as before.

\begin{defn}
A connected component $V$ of $\tilde{\cU}_{\II'}$
is called a {\it $\II'$-region} of $(\gamma; \delta)$
if it satisf\/ies the following conditions:
\begin{enumerate}[(i)]\itemsep=0pt
\item $V$ contains at least one $\II'$-unit of height $0$.
\item $V$ does not contain any boundary $\II'$-unit.
\end{enumerate}
\end{defn}

In Figs.~\ref{fig:extra-rule-example} and \ref{fig:gamma},
we can see that the $\II$-region of $(\alpha; \beta)$ and
the $\II'$-region of $(\gamma; \delta)$ coincide.
In fact,

\begin{prop}\label{prop:equivalence-II}
Let $(\alpha;\beta)=\pi(p)$ and $(\gamma;\delta)=\pi'(\bp)$.
Then,
$V$ is a $\II$-region of $(\alpha;\beta)$ if and only
if
$V$ is a $\II'$-region of $(\gamma;\delta)$.
\end{prop}

\begin{proof}
It is not dif\/f\/icult to show the following facts:
\begin{itemize}\itemsep=0pt
\item[(1)]
A $\II'$-unit of $(\gamma;\delta)$ is a $\II$-
unit of
$(\alpha;\beta)$.
\item[(2)] Let $U$ be a $\II$-unit of $(\alpha;
\beta)$
satisfying one of the following conditions:
\begin{itemize}\itemsep=0pt
\item[(i)] $U\subset S_+$, and  \eqref{eq:II-unit}
holds for some $i$ with $r<i<l$.
\item[(ii)] $U\subset S_-$, and  \eqref{eq:II-unit}
holds for some $i$ with $0<i<l-r$.
\end{itemize}
Then, $U$ is a $\II'$-unit of $(\gamma;\beta)$.
\item[(3)]
Let $U$ be a $\II$-unit of $(\alpha;\beta)$
satisfying one of the following conditions:
\begin{itemize}\itemsep=0pt
\item[(i)] $U\subset S_+$, and  \eqref{eq:II-unit}
holds for some $i$ with $i\leq r$ or $i=l$.
\item[(ii)] $U\subset S_-$, and  \eqref{eq:II-unit}
holds for some $i$ with $i=0$ or $i\geq l-r$.
\end{itemize}
Let $U'$ be a unit which belongs to a $\II'$-region of
$(\gamma;\delta)$.
Then, $U$ and $U'$ are not adjacent to each other.
\end{itemize}
Now, the if part of the proposition follows from
(1), (2), and (3), while the only if part
follows from (1) and (2).
\end{proof}

\subsection{Extra rule in terms of tableau}
With Proposition \ref{prop:equivalence-II},
it is now straightforward to translate the extra
rule $(\bE)$
into tableau language as type $D_n$ \cite{NN2}.
We only give the result.

Fix an HV-tableau $T$ of shape $\lambda/\mu$.
For any $a_1, \dots, a_m \in I$,
let $C(a_1, \dots, a_m)$ be a conf\/i\-gu\-ra\-tion
in $T$ as follows:
\begin{equation}
\text{
{\setlength{\unitlength}{0.2mm}
\begin{picture}(25,100)
\multiput(0,0)(25,0){2}{\line(0,1){110}}
\multiput(0,0)(0,25){2}{\line(1,0){25}}
\multiput(0,110)(0,-25){3}{\line(1,0){25}}
\put(4,92){$a_1$}
\put(4,67){$a_2$}
\put(9,34){$\vdots$}
\put(2,7){$a_{\scriptscriptstyle m}$}
\end{picture}
}
}
\end{equation}
We call it an {\it L-configuration} if it satisf\/ies
\begin{itemize}
\item[(i)]
$1 \preceq a_1 \prec \dots \prec a_m \preceq n$.
\item[(ii)]
If $a_m=n$, $\Pos(a_m)=(i,j)$, and $(i,j-1)\in\lambda/\mu$,
then $T(i,j-1)\neq \overline{n}$.
\end{itemize}
Here and below, $\Pos(a_m)$, for example, means the position of $a_m$ in $T$.
We call it a {\it U-configuration} if it satisf\/ies
\begin{itemize}
\item[(i)]
$\overline{n} \preceq a_1 \prec \dots \prec a_m \preceq \overline{1}$.
\item[(ii)]
If $a_1=\overline{n}$, $\Pos(a_1)=(i,j)$,
and $(i,j+1)\in\lambda/\mu$,
then $T(i,j+1)\neq n$.
\end{itemize}
An L-conf\/iguration is identif\/ied with a part of lower paths
for $T$  under the map $\cTh$, while a U-conf\/iguration
is so with  a part of upper paths.

\begin{figure}[t]
\centerline{\includegraphics{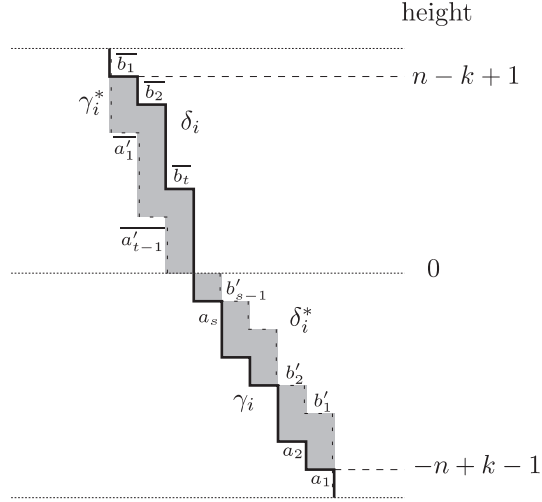}}
\caption{An example of a pair of lower and upper $e$-paths
$(\gamma_i,\delta_i)$ whose
part corresponds to a LU-conf\/iguration
of type 1 as in \eqref{eq:type-one}.}\label{fig:tableau-path-correspondence}
\end{figure}

Let $(L,U)$ be
a pair of an L-conf\/iguration $L=C(a_1, \dots, a_s)$
and a U-conf\/iguration
$U=C(\overline{b_t}, \dots, \overline{b_1})$ in the same column.
We call it an {\it LU-configuration} of $T$ if it satisf\/ies
one of the following two conditions:

Condition 1. {\it LU-configuration of type 1}.
$(L,U)$ has the form
\begin{equation}\label{eq:type-one}
 \raisebox{-50pt}
{{\setlength{\unitlength}{0.2mm}
\begin{picture}(50,180)
{\thicklines
\put(0,0){\line(0,1){180}}
\put(25,0){\line(0,1){180}}
\put(0,0){\line(1,0){25}}
\put(0,25){\line(1,0){25}}
\put(0,65){\line(1,0){25}}
\put(0,90){\line(1,0){25}}
\put(0,115){\line(1,0){25}}
\put(0,155){\line(1,0){25}}
\put(0,180){\line(1,0){25}}
}
\put(4,165){$a_1$}
\put(10,127){$\vdots$}
\put(4,100){$a_s$}
\put(4,5){$\overline{b_1}$}
\put(10,37){$\vdots$}
\put(4,70){$\overline{b_t}$}
\multiput(25,0)(5,0){15}{\line(1,0){2}}
\multiput(25,180)(5,0){15}{\line(1,0){2}}
\put(80,0){\vector(0,1){180}}
\put(80,180){\vector(0,-1){180}}
\put(85,100){$\scriptstyle n-k+2$}
\end{picture}
}
}
\end{equation}
for some $k$ with $1 \le k \le n$,
$n-k+2=s+t$, and
\begin{gather}
  a_1 =  k, \qquad
  \overline{b_1} = \overline{k}, \\
a_{i+1} \preceq b'_i, \quad (1 \le i \le s-1), \qquad
\overline{b_{i+1}} \succeq \overline{a'_i}, \quad (1 \le i \le t-1),
\label{eq:intersecting-cond}
\end{gather}
where $a'_1 \prec \dots \prec a'_{t}$
and $b'_1 \prec \dots \prec b'_{s}$
are def\/ined as
\begin{gather}
\{a_1, \dots, a_s \} \sqcup \{ a'_1, \dots, a'_{t-1} \}
 = \{ k, k+1, \dots, n \}, \quad a'_t=\overline{n}, \nonumber\\
\{ b_1, \dots, b_{t} \} \sqcup
\{ b'_1, \dots, b'_{s-1} \}
 = \{ k, k+1, \dots, n \},
\quad b'_s = \overline{n}.\label{eq:dual-type-one}
\end{gather}
See Fig.~\ref{fig:tableau-path-correspondence}
for the corresponding part in the paths.

\begin{figure}
\centerline{\includegraphics{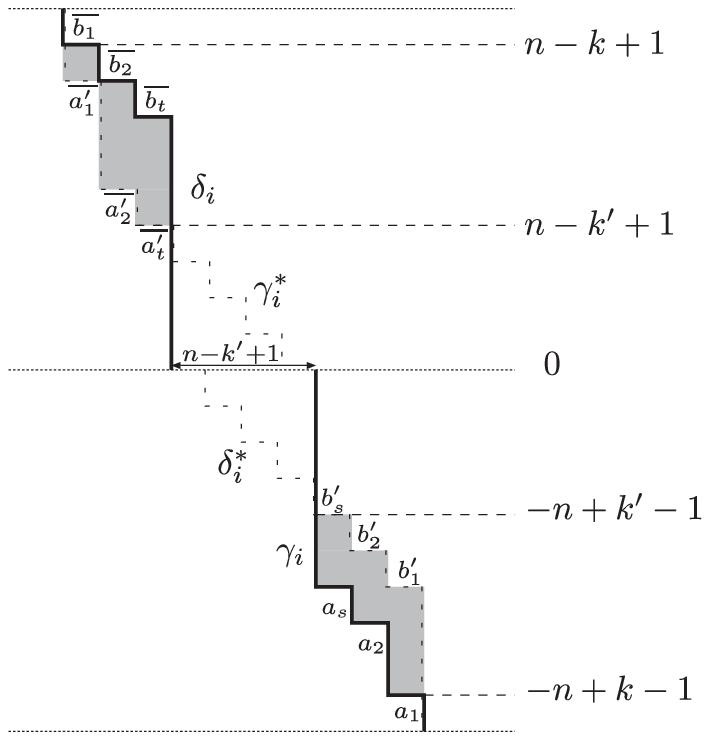}}
\caption{An example of a pair of lower and upper $e$-paths
$(\gamma_i, \delta_i)$ whose
part corresponds to a LU-conf\/iguration
of type 2 as in \eqref{eq:type-two}.}
\label{fig:tableau-path-correspondence-two}
\end{figure}

Condition 2. {\it LU-configuration of type 2}.
$(L,U)$ has the form
\begin{equation}\label{eq:type-two}
\hspace*{18mm} \raisebox{-58pt}{
{\setlength{\unitlength}{0.2mm}
\begin{picture}(50,210)
{\thicklines
\put(0,130){\line(0,1){80}}
\put(25,130){\line(0,1){80}}
\put(0,0){\line(0,1){80}}
\put(25,0){\line(0,1){80}}
\put(0,130){\line(1,0){25}}
\put(0,155){\line(1,0){25}}
\put(0,185){\line(1,0){25}}
\put(0,210){\line(1,0){25}}
\put(0,0){\line(1,0){25}}
\put(0,25){\line(1,0){25}}
\put(0,55){\line(1,0){25}}
\put(0,80){\line(1,0){25}}
}
\put(5,195){$a_1$}
\put(10,160){$\vdots$}
\put(5,140){$a_s$}
\put(8,115){$a$}
\put(6,4){$\overline{b_1}$}
\put(10,30){$\vdots$}
\put(6,59){$\overline{b_t}$}
\put(9,82){$b$}
\multiput(25,0)(5,0){15}{\line(1,0){2}}
\multiput(25,210)(5,0){15}{\line(1,0){2}}
\put(80,0){\vector(0,1){210}}
\put(80,210){\vector(0,-1){210}}
\put(85,90){$\scriptstyle n-k+2$}
\multiput(0,80)(-5,0){10}{\line(-1,0){2}}
\multiput(0,130)(-5,0){10}{\line(-1,0){2}}
\put(-30,80){\vector(0,1){50}}
\put(-30,130){\vector(0,-1){50}}
\put(-100,100){$\scriptstyle n-k'+1$}
%
%\multiput(0,80)(0,5){10}{\line(0,1){2}}
%\multiput(25,80)(0,5){10}{\line(0,1){2}}
\end{picture}
}
}
\end{equation}
for some $k$ and $k'$ with $1 \le k  < k' \le n$,
$n-k+2=n-k'+1+s+t$, and
\begin{gather}
a_1 =k, \qquad
\overline{b_1} = \overline{k}, \qquad
a'_{t}=k', \qquad
\overline{b'_{s}}=\overline{k'},\\
a \succ k' \qquad  \mbox{or}\qquad
\bigg[
\begin{gathered}
\mbox{if}\ a=n, \Pos(a)=(i,j), \ \mbox{then}\ \\
(i,j-1)\in\lambda/\mu, T(i,j-1)=\overline{n}
\end{gathered}
\bigg],
\\
b \prec \overline{k'} \qquad  \mbox{or}\qquad
\bigg[
\begin{gathered}
\mbox{if}\
b=\overline{n}, \Pos(b)=(i,j), \ \mbox{then}\
\\
(i,j+1)\in\lambda/\mu, T(i,j+1)=n
\end{gathered}
\bigg],
\\
a_{i+1} \preceq b'_i, \quad (1 \le i \le s-1), \qquad
\overline{b_{i+1}} \succeq \overline{a'_i}, \quad
(1 \le i \le t-1),
\end{gather}
where
$a'_1 \prec \dots \prec a'_{t}$  and
$b'_1 \prec \dots \prec b'_{s}$  are def\/ined as
\begin{gather}
\{a_1, \dots, a_s \} \sqcup \{ a'_1, \dots, a'_{t} \}
 = \{ k, k+1, \dots, k'\}, \nonumber\\
\{ b_1, \dots, b_{t} \} \sqcup
\{ b'_1, \dots, b'_{s} \}
 = \{ k, k+1, \dots,  k' \}.\label{eq:dual-type-two}
\end{gather}
See Fig.~\ref{fig:tableau-path-correspondence-two}
for the corresponding parts in the paths.

We say that an $L$-conf\/iguration $L=C(a_1, \dots, a_m)$
in the $j$th column of $T$ is {\it boundary} if
$\Pos(a_1)=(\mu'_j+1,j)$, i.e., if $a_1$ is at the top of the
$j$th column, and $m$ is the largest number such that
$L \cap L'=\varnothing$ for any $LU$-conf\/iguration $(L',U')$.
Similarly, a $U$-conf\/iguration $U=C(a_1, \dots, a_m)$
in the $j$th column of $T$ is {\it boundary} if
$\Pos(a_m)=(\lambda'_j,j)$, i.e., if $a_m$ is at the bottom of the
$j$th column, and $m$ is the largest number such that
$U \cap U'=\varnothing$ for any $LU$-conf\/iguration $(L',U')$.

Let $(L,U)=(C(a_1, \dots, a_s),
C(\overline{b_t}, \dots, \overline{b_1}))$ be
an LU-conf\/iguration, and set
$a'_1\prec \dots \prec a'_{t}$ and
$b'_1\prec \dots \prec b'_{s}$ as in
\eqref{eq:dual-type-one} (resp.~as in \eqref{eq:dual-type-two})
if $(L,U)$ is of type 1 (resp.~of type 2).
We say that an L-conf\/iguration
$L'$ is {\it right-adjacent} to $(L,U)$ if
$L'$ is in the right-next column to~$L$; furthermore,
there exists some pair of an entry $e$ of $L'$ and an entry $a_i$ of~$L$
such that $e$ is right-next to $a_i$ and $e\prec b'_i$.
Similarly,
we say that a U-conf\/iguration $U'$
is {\it left-adjacent} to $(L,U)$ if
$U'$ is in the left-next column to $U$; furthermore,
there exists some pair of an entry $e$ of $U'$ and an entry $\overline{b_i}$
of~$U$ such that $e$ is left-next to $\overline{b_i}$
and $e\succ \overline{a'_i}$,
where $\overline{\overline{n}}=n$.
Then, we say that an LU-conf\/iguration
$(L',U')$ is {\it adjacent} to $(L,U)$ if one of the following
is satisf\/ied, and write it by $(L, U)\diamond (L',U')$:
\begin{enumerate}[(i)]\itemsep=0pt
\item $L'$ is right-adjacent to $(L,U)$.
\item $L$ is right-adjacent to $(L',U')$.
\item $U'$ is left-adjacent to $(L,U)$.
\item $U$ is left-adjacent to $(L',U')$.
\end{enumerate}

For any tableau $T$, let $\cLU(T)$ be
the set of all LU-conf\/igurations of $T$.
Then, the adjacent relation $\diamond$ of the LU-conf\/igurations
generates an equivalence relation $\sim$ in $\cLU(T)$.

\begin{defn}
For any $(L,U) \in \cLU(T)$,
let $[(L,U)]\subset \cLU(T)$ be the equivalence class of $(L, U)$
with respect to $\sim$,
and let $R=R(L,U):=\bigcup_{(L',U') \in [(L,U)]}(L',U')$
be the corresponding conf\/iguration in $T$.
We call $R$  a {\it $\II$-region} of $T$, if
the following is satisf\/ied:
\begin{enumerate}\itemsep=0pt
\item[(i)]
No boundary L-conf\/iguration $L$
is right-adjacent to $L'$ for any
LU-conf\/iguration $(L',U')$ in $R$.
\item[(ii)]
No boundary U-conf\/iguration
$U$ is left-adjacent to $U'$ for any
LU-conf\/iguration $(L',U')$ in~$R$.
\end{enumerate}
Moreover,
we say that $R$ is {\it odd\/} if
the number of the type 1 LU-conf\/igurations
in $R$ is odd.
\end{defn}

Then, an odd $\II$-region of $T=\cTh(\bp)$
corresponds to an odd $\II$-region of $\bp$, and therefore,
Theorem \ref{thm:tableau-description} is rewritten as follows:

\begin{thm}[Tableaux description, cf.\ {\cite[Theorem 5.6]{NN2}}]
\label{thm:configuration}
For any skew diagram $\lambda/\mu$ satisfying
the positivity condition \eqref{eq:positivity},
we have
\begin{equation*}
\chi_{\lambda/\mu, a} = \sum_{T \in \Tab(\lambda/\mu)} z_a^{T},
\end{equation*}
where $\Tab(\lambda/\mu)$ is the set of all
the tableaux of shape $\lambda/\mu$ which satisfy the
$(\bH)$, $(\bV)$,
and the following extra rule $(\bE')$:
\vspace*{5pt}

\noindent
\begin{tabular}{ll}
$(\bE')$ & $T$ does not have any odd $\II$-region.
\end{tabular}
\vspace{5pt}
\end{thm}

\begin{exmp} Consider the tableaux $T$ and $T'$ in Example
\ref{exmp:two-tableaux}.
The conf\/iguration
\begin{equation}
\raisebox{-26pt}
{
{\setlength{\unitlength}{0.2mm}
\begin{picture}(25,100)
\multiput(25,0)(25,0){3}{\line(0,1){100}}
\multiput(0,0)(0,25){4}{\line(1,0){75}}
\put(0,0){\line(0,1){75}}
\put(25,100){\line(1,0){50}}
%\multiput(0,100)(0,-25){3}{\line(1,0){25}}
%\put(7,80){$1$}
\put(7,55){$3$}
\put(7,30){$\overline{4}$}
\put(7,5){$\overline{3}$}
\put(32,80){$2$}
\put(32,55){$4$}
\put(32,30){$\overline{4}$}
\put(32,5){$\overline{2}$}
\put(57,80){$2$}
\put(57,55){$4$}
\put(57,30){$\overline{3}$}
\put(57,5){$\overline{2}$}
\end{picture}
}
}
\hskip25pt
\subset T
\end{equation}
is an odd $\II$-region of $T$, while the conf\/iguration
(the boxed part in the following)
\begin{equation}
\raisebox{-26pt}
{
{\setlength{\unitlength}{0.2mm}
\begin{picture}(25,100)
\multiput(50,0)(25,0){2}{\line(0,1){100}}
\multiput(0,25)(0,25){3}{\line(1,0){75}}
\put(0,0){\line(1,0){25}}
\put(50,0){\line(1,0){25}}
\put(0,0){\line(0,1){75}}
\put(25,0){\line(0,1){75}}
\put(50,100){\line(1,0){25}}
%\multiput(0,100)(0,-25){3}{\line(1,0){25}}
%\put(7,80){$1$}
\put(7,55){$3$}
\put(7,30){$\overline{4}$}
\put(7,5){$\overline{3}$}
%\put(32,80){$2$}
\put(32,55){$4$}
\put(32,30){$\overline{4}$}
\put(32,5){$\overline{2}$}
\put(57,80){$2$}
\put(57,55){$4$}
\put(57,30){$\overline{3}$}
\put(57,5){$\overline{2}$}
\end{picture}
}
}
\hskip25pt
\subset T'
\end{equation}
is not a $\II$-region of $T'$ because
the boundary $U$-conf\/iguration $C(\overline{2})$
in the second column is
left-adjacent to the LU-conf\/iguration
$(C(2,4), C(\overline{3},\overline{2}))$
in the third column.
\end{exmp}

\begin{exmp}
Let $\lambda/\mu$ be a skew diagram of at most
three rows.
The following is
an example of an odd $\II$-region of $T\in \HVTab(\lambda/\mu)$:
\begin{equation}
\raisebox{-25pt}{
\begin{picture}(60,60)
\put(0,0){\line(1,0){60}}
\put(0,20){\line(1,0){60}}
\put(0,40){\line(1,0){60}}
\put(20,60){\line(1,0){40}}
\put(0,0){\line(0,1){40}}
\put(20,0){\line(0,1){60}}
\put(40,0){\line(0,1){60}}
\put(60,0){\line(0,1){60}}
\put(67,47){$ a$}
\put(-13,7){$ b$}
\put(7,27){$ n$}
\put(7,7){$ \overline{n}$}
\put(23,47){$ \scriptstyle n-1$}
\put(27,27){$ n$}
\put(23,7){$ \scriptstyle \overline{n-1}$}
\put(43,47){$ \scriptstyle n-1$}
\put(47,27){$ \overline{n}$}
\put(43,7){$ \scriptstyle \overline{n-1}$}
\end{picture}
}
\hskip10pt,
\end{equation}
where $a\succeq \overline{n}$ and $b\preceq n$
if they exist.
The complete list of all the possible odd $\II$-regions
in $T$ corresponds to the rules ({\bf E-2R}) and ({\bf E-3R})
of Theorem 5.7 in \cite{NN1}.
\end{exmp}

\begin{exmp}
Let $\lambda/\mu$ be a skew diagram of two columns
satisfying the positivity condition~\eqref{eq:positivity}.
We note that if $T\in \HVTab(\lambda/\mu)$
contains a type 1 LU-conf\/iguration, say,
in the f\/irst column, then $\gamma_2(0)=\delta_2(0)$.
Then, it is easy to prove that the extra rule
$(\bE')$ is equivalent
to the following condition:

\smallskip
({\bf E-2C}) \quad $T$ does not have any
odd type 1 LU-conf\/iguration as
\begin{equation*}
\raisebox{-50pt}
{
{\setlength{\unitlength}{0.2mm}
\begin{picture}(50,180)
{\thicklines
\put(0,0){\line(0,1){180}}
\put(25,0){\line(0,1){180}}
\put(0,0){\line(1,0){25}}
\put(0,25){\line(1,0){25}}
\put(0,65){\line(1,0){25}}
\put(0,90){\line(1,0){25}}
\put(0,115){\line(1,0){25}}
\put(0,155){\line(1,0){25}}
\put(0,180){\line(1,0){25}}
}
\put(4,163){$a_1$}
\put(10,127){$\vdots$}
\put(4,98){$a_s$}
\put(4,5){$\overline{b_1}$}
\put(10,37){$\vdots$}
\put(4,70){$\overline{b_t}$}
\put(29,163){$d_1$}
\put(35,127){$\vdots$}
\put(29,98){$d_s$}
\put(-21,5){$c_1$}
\put(-15,37){$\vdots$}
\put(-21,70){$c_t$}
%
%\multiput(25,0)(5,0){15}{\line(1,0){2}}
%\multiput(25,180)(5,0){15}{\line(1,0){2}}
%\put(80,0){\vector(0,1){180}}
%\put(80,180){\vector(0,-1){180}}
%\put(85,100){$\scriptstyle n-k+2$}
%
\end{picture}
}
},
\end{equation*}
where
\begin{enumerate}[(i)]\itemsep=0pt
\item $(L,U)=(C(a_1, \dots, a_s),
C(\overline{b_t}, \dots, \overline{b_1}))$
is a type 1 LU-conf\/iguration.
\item Let $a'_i$ be the one in \eqref{eq:dual-type-one}.
For each $k=1, \dots, t$, if $\Pos(\overline{b_k})=(i,j)$ and
$(i,j-1) \in \lambda/\mu$, then
$c_k:=T(i,j-1) \preceq \overline{a'_k}$.
\item Let $b'_i$ be the one in \eqref{eq:dual-type-one}.
For each $k=1, \dots, s$, if $\Pos(a_k)=(i,j)$ and
$(i,j+1) \in \lambda/\mu$, then
$d_k:=T(i,j+1) \succeq b'_k$.
\end{enumerate}
\smallskip

This proves Conjecture 5.9 in \cite{NN1}.
\end{exmp}

\appendix

\section{Remarks on Section \ref{subsec:positive-sum}}\label{sec:appendix}

We give some technical remarks on
Section \ref{subsec:positive-sum}.

\subsection{Expansion and folding}

The involution $\iota_2$ in Proposition \ref{prop:second-involution}
is def\/ined by using the deformations of paths called
{\it expansion} and {\it folding}.
They are generalizations of the `resolution of a (\ref{eq:positive-sum})
pair of paths' in \cite{NN1} and its inverse.
They are def\/ined  in exactly the same way
as type $D_n$ \cite{NN2}.

For any $(\alpha; \beta)\in \cH(\lambda/\mu)$,
let $V$ be any $\I$- or $\II$-region of $(\alpha; \beta)$.
Let $\alpha'_i$ be the lower path obtained from
$\alpha_i$ by replacing the part $\alpha_i\cap V$ with $\beta^*_{i+1}\cap V$,
and let $\beta'_i$ be the upper path obtained from
$\beta_i$ by replacing the part $\beta_i\cap V$ with $\alpha^*_{i-1}\cap V$.
Set
$\varepsilon_V(\alpha; \beta):=
(\alpha'_1, \dots, \alpha'_l; \beta'_1, \dots, \beta'_l)$.
We have
\begin{prop}\label{prop:I-II-region}
Let $\lambda/\mu$ be a skew diagram
satisfying the positivity condition \eqref{eq:positivity}.
Then, for any $(\alpha; \beta) \in \cH(\lambda/\mu)$,
we have
\begin{enumerate}\itemsep=0pt
\item \label{item:prop-one}
For any $\I$- or $\II$-region
$V$ of $(\alpha; \beta)$,
$\varepsilon_V(\alpha; \beta) \in \cH(\lambda/\mu)$.
\item \label{item:prop-two}
For any $\I$-region $V$ of $(\alpha; \beta)$,
$V$ is a $\II$-region of $\varepsilon_V(\alpha; \beta)$.
\item For any $\II$-region $V$ of $(\alpha; \beta)$,
$V$ is a $\I$-region of $\varepsilon_V(\alpha; \beta)$.
\end{enumerate}
\end{prop}

\begin{figure}
\centerline{\includegraphics{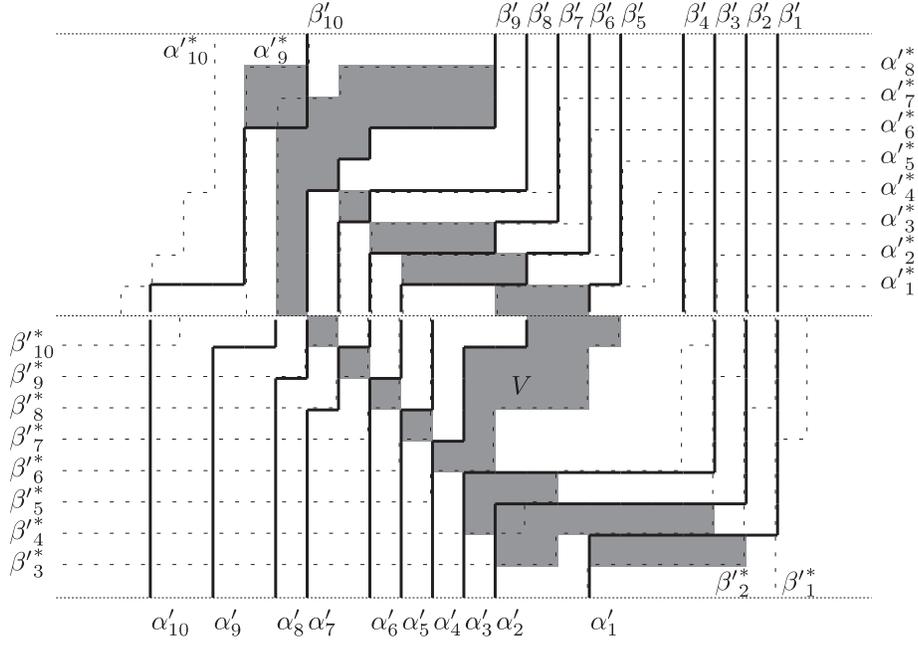}}
\caption{$(\alpha'; \beta')$ is the folding of $(\alpha; \beta)$
in Figure \ref{fig:I-region} with respect to $V$.
Conversely, $(\alpha; \beta)$ is the expansion of
$(\alpha'; \beta')$ with respect to $V$.}\label{fig:expansion}
\end{figure}

We call the correspondence $(\alpha; \beta) \mapsto \varepsilon_V(\alpha; \beta)$
the {\it expansion} (resp.~the {\it folding}) with respect to $V$,
if $V$ is a $\I$-region (resp.~a $\II$-region) of $(\alpha; \beta)$.
See Fig.~\ref{fig:expansion} for an example.

\subsection[$\I_k$- and $\II_k$-units]{$\boldsymbol{\I_k}$- and $\boldsymbol{\II_k}$-units}\label{sec:folding-map}

To construct the folding map $\phi$
in Proposition \ref{prop:folding-map},
which is a key to derive the tableaux description,
we generalize the expansion and the folding
to the {\it $k$-expansion} and the {\it $k$-folding} \cite{NN2}.
The original corresponds to $k=1$.
Like the $k=1$ case, it starts from the def\/initions
of $\I_k$- and $\II_k$-units, which are slightly modif\/ied for type $C_n$.
\begin{defn}\label{def:k-unit}
Let
$(\alpha; \beta) \in \cH(\lambda/\mu)$.
For any unit $U \subset S_{\pm}$, let $\pm r= \Ht(U)$ and
let $a$ and $a'=a+1$ be the horizontal position of
the left and right edges of $U$.
Then, for any $k=1 ,\dots , l-1$,
\begin{enumerate}\itemsep=0pt
\item
$U$ is called a $\I_k$-{\it unit} of $(\alpha; \beta)$
if there exists some $i$ ($0\le i \le l$) such that
\begin{gather}
 \alpha^*_i(r) \le a < a' \le \beta_{i+k}(r), \quad \text{if $U \subset S_+$}, \nonumber\\
 \alpha_i(-r) \le a < a' \le \beta_{i+k}^*(-r), \quad \text{if $U \subset S_-$}.\label{eq:Ik-unit}
\end{gather}
\item
$U$ is called a $\II_k$-{\it unit} of $(\alpha; \beta)$
if there exists some $i$ ($0 \le i \le l$) such that
\begin{gather}
 \beta_{i+k}(r) \le a < a' \le \alpha^*_i(r), \quad \text{if $U \subset S_+$}, \nonumber\\
\beta^*_{i+k}(-r) \le a < a' \le \alpha_i(-r), \quad \text{if $U \subset S_-$}.\label{eq:IIk-unit}
\end{gather}
\end{enumerate}
Here, we set $\beta_i(r)=\beta^*_i(-r)=-\infty$
and $\alpha_i(-r)=\alpha^*_i(r)=+\infty$  for any $r$
and for any $i \ne 1, \dots, l$.
Furthermore, a $\II$-unit $U$ of $(\alpha; \beta)$
is called a {\it boundary} $\II$-unit if
one of the following holds:
\begin{itemize}\itemsep=0pt
\item[(i)] $U\subset S_+$, and  \eqref{eq:IIk-unit}
holds for $i\leq r $, $i\ge l+1-k$, or $r=n$.
\item[(ii)] $U\subset S_-$, and  \eqref{eq:IIk-unit}
holds for $i=0 $, $i\geq l+1-k-r$, or $r=n$.
\end{itemize}
\end{defn}
Then, the rest of the def\/initions and
the proof of Proposition \ref{prop:folding-map}
are exactly the same as type~$D_n$.

\pdfbookmark[1]{References}{ref}
\LastPageEnding

\end{document}